\documentclass[11pt, a4paper]{amsart}
\usepackage{mathbbol}
\usepackage{amsmath, amsthm, amssymb}
\usepackage{txfonts, mathrsfs}
\usepackage[usenames]{color}
\usepackage{psfrag}
\usepackage{graphicx}
\usepackage{todonotes}

\newtheorem*{thm*}{Theorem}

\newtheorem{thm}{Theorem}[section]
\newcommand{\bt}{\begin{thm}}
\newcommand{\et}{\end{thm}}

\newtheorem{cor}[thm]{Corollary}
\newcommand{\bc}{\begin{cor}}
\newcommand{\ec}{\end{cor}}

\newtheorem{lem}[thm]{Lemma}
\newcommand{\bl}{\begin{lem}}
\newcommand{\el}{\end{lem}}

\newtheorem{prop}[thm]{Proposition}
\newcommand{\bp}{\begin{prop}}
\newcommand{\ep}{\end{prop}}

\newtheorem{defn}[thm]{Definition}
\newcommand{\bd}{\begin{defn}}      
\newcommand{\ed}{\end{defn}}

\newtheorem{rmrk}[thm]{Remark}
\newcommand{\br}{\begin{rmrk}}
\newcommand{\er}{\end{rmrk}}

\newtheorem{quest}[thm]{Question}
\newcommand{\bq}{\begin{quest}}
\newcommand{\eq}{\end{quest}}

\newtheorem{example}[thm]{Example}

\newcommand{\N}{\mathbb{N}}

\newcommand{\R}{\mathbb{R}}

\newdimen\vintkern\vintkern12pt
\def\vint{-\kern-\vintkern\int}

\newcommand{\hm}{{\mathcal H}}

\newcommand{\diam}{\operatorname{diam}}
\newcommand{\trace}{\operatorname{tr}}
\newcommand{\length}{\ell}
\newcommand{\Area}{\operatorname{Area}}

\newcommand{\osc}{\operatorname{osc}}

\newcommand{\md}{\operatorname{md}}

\newcommand{\jac}{{\mathbf J}}
\newcommand{\ap}{\operatorname{ap}}
\newcommand{\apmd}{\md}

\DeclareMathOperator{\interior}{int}

\DeclareMathOperator{\MOD}{mod}

\begin{document}
\bibliographystyle{plain}

\title[Quasiconformal almost parametrizations]{Quasiconformal almost parametrizations of metric surfaces}

\author{Damaris Meier}


\address
  {Department of Mathematics\\ University of Fribourg\\ Chemin du Mus\'ee 23\\ 1700 Fribourg, Switzerland}
\email{damaris.meier@unifr.ch}

\author{Stefan Wenger}

\address
  {Department of Mathematics\\ University of Fribourg\\ Chemin du Mus\'ee 23\\ 1700 Fribourg, Switzerland}
\email{stefan.wenger@unifr.ch}

\date{\today}

\thanks{Research supported by Swiss National Science Foundation Grant 182423.}

\begin{abstract}
 We look for minimal conditions on a two-dimensional metric surface $X$ of locally finite Hausdorff $2$--measure under which $X$ admits an (almost) parametrization with good geometric and analytic properties. Only assuming that $X$ is locally geodesic, we show that Jordan domains in $X$ of finite boundary length admit a quasiconformal almost parametrization. If $X$ satisfies some further conditions then such an almost parametrization can be upgraded to a geometrically quasiconformal homeomorphism or a quasisymmetric homeomorphism. In particular, we recover Rajala's recent quasiconformal uniformization theorem in the special case that $X$ is locally geodesic as well as Bonk-Kleiner's quasisymmetric uniformization theorem. On the way we establish the existence of Sobolev discs spanning a given Jordan curve in $X$ under nearly minimal assumptions on $X$ and prove the continuity of energy minimizers.
\end{abstract}

\maketitle

\section{Introduction and statement of main results}

\subsection{Background}

Every smooth Riemann surface is conformally diffeomorphic to a surface of constant curvature by the classical uniformization theorem. The uniformization problem for metric spaces, widely studied in the field of analysis in metric spaces and of importance also in other areas, asks to find conditions on a given metric space $X$, homeomorphic to some model space $M$, under which there still exists a homeomorphism from $X$ to $M$ with good geometric and analytic properties. 

In this paper we consider the uniformization problem for metric spaces homeomorphic to a two-dimensional surface and of locally finite Hausdorff $2$--measure. In this setting, two outstanding uniformization results were proved in \cite{BK02} and \cite{Raj14}. Bonk-Kleiner \cite{BK02} showed that an Ahlfors $2$--regular metric space $X$ homeomorphic to the standard two-sphere $S^2$ admits a quasisymmetric homeomorphism to $S^2$ if and only if $X$ is linearly locally connected. Ahlfors $2$--regular means that the Hausdorff $2$--measure of balls of radius $r$ is comparable to $r^2$ and a quasisymmetric homeomorphism is a homeomorphism that distorts shapes in a controlled manner. We refer to Section~\ref{sec:properties-qc-almost-parametrizations} for the definitions of quasisymmetric homeomorphism and linear local connectedness. More recently, Rajala \cite{Raj14} gave a characterization of metric planes admitting a geometrically quasiconformal homeomorphism to a Euclidean domain. His characterization involves a condition called reciprocality. A geometrically quasiconformal map is a homeomorphism that leaves the (conformal) modulus of curve families invariant up to a multiplicative constant.
Rajala's result in particular gives a new approach to the Bonk-Kleiner quasisymmetric uniformization theorem. The results in \cite{BK02} and \cite{Raj14} have been extended for example in \cite{BK05}, \cite{Wil08}, \cite{Wil10}, \cite{MW13}, \cite{Iko19}. In \cite{LW20-param}, Lytchak and the second author provided a further approach to the Bonk-Kleiner theorem which relies on results about the existence and regularity of energy and area minimizing discs in metric spaces admitting a quadratic isoperimetric inequality developed in \cite{LW15-Plateau} and \cite{LW16-intrinsic}.

While we work with metric surfaces of locally finite Hausdorff $2$--measure in this paper, the uniformization problem has also been studied for fractal spaces, see for example \cite{Mey02}, \cite{Mey10}, \cite{LRR17}, \cite{RRR19}, \cite{RRasim20}.
The aim of our paper is to establish the existence of parametrizations or almost parametrizations with good properties under nearly minimal conditions on $X$. The properties are such that they upgrade to geometrically quasiconformal parametrizations under Rajala's reciprocality condition and to quasisymmetric parametrizations under the condition of Ahlfors $2$--regularity and linear local connectedness. On the way to prove our parametrization results, we establish the existence of Sobolev discs spanning a given Jordan curve under nearly minimal assumptions on $X$ and regularity of energy minimizers.

\subsection{Parametrization results}

We now turn to a rigorous discussion of our results. Let $X$ be a metric space homeomorphic to a two-dimensional surface and assume that $X$ has locally finite Hausdorff $2$--measure. The modulus of a family $\Gamma$ of curves in $X$ is defined by $$\MOD(\Gamma):= \inf_\rho \int_X\rho^2\,d\hm^2,$$ where the infimum is taken over all Borel functions $\rho\colon X\to [0,\infty]$ for which $\int_\gamma\rho\geq 1$ for every $\gamma\in\Gamma$, see Section~\ref{sec:prelims}. A homeomorphism $u\colon D\to \Omega$ from the unit disc $D\subset\R^2$ to a domain $\Omega\subset X$ is called geometrically quasiconformal if $u$ leaves the modulus of curve families quasi-invariant, thus there exists $K\geq 1$ such that 
\begin{equation}\label{eq:geom-qc-def-intro}
 K^{-1}\cdot \MOD(\Gamma)\leq \MOD(u\circ\Gamma)\leq K\cdot \MOD(\Gamma)
\end{equation}
for every family $\Gamma$ of curves in $D$. Here, $u\circ\Gamma$ denotes the family of curves $u\circ\gamma$ with $\gamma\in\Gamma$. By Rajala's recent uniformization result \cite{Raj14}, $\Omega$ admits a geometrically quasiconformal parametrization if and only if $\Omega$ satisfies a certain reciprocality condition described below. It is natural to wonder to what extent one can weaken this condition and still obtain an (almost) parametrization of $\Omega$ with suitable properties which can then be upgraded to a geometrically quasiconformal homeomorphism when $X$ is reciprocal. Our main result shows that the reciprocality condition can be dropped completely, at least when the underlying metric space $X$ is locally geodesic.

\bt\label{thm:main-almost-param}
 Let $X$ be a locally geodesic metric space homeomorphic to $\R^2$ and of locally finite Hausdorff $2$--measure. If $\Omega\subset X$ is a Jordan domain of finite boundary length then there exists a continuous, monotone surjection $u\colon \overline{D}\to \overline{\Omega}$ such that 
 \begin{equation}\label{eq:analytic-qc-mod}
  \MOD(\Gamma) \leq K\cdot\MOD(u\circ \Gamma)
 \end{equation}
 for every family $\Gamma$ of curves in $\overline{D}$, where $K=\frac{4}{\pi}$.
\et

The map $u$ is called monotone if $u^{-1}(x)$ is connected for every $x\in X$; equivalently, $u$ is the uniform limit of homeomorphisms $u_n\colon \overline{D}\to \overline{\Omega}$. Notice that in the generality of Theorem~\ref{thm:main-almost-param}, there need not exist a homeomorphism satisfying \eqref{eq:analytic-qc-mod}, see Example~\ref{ex:collapsed-disc}. When $u$ is a homeomorphism then \eqref{eq:analytic-qc-mod} is equivalent to the so-called analytic definition of quasiconformality by \cite{Wil12}. Theorem~\ref{thm:main-almost-param} answers a question of Rajala and the second author stated e.g.~in \cite[Question 1.1]{IR20}
in the special case of locally geodesic metric spaces. The factor $K=\frac{4}{\pi}$ appearing in Theorem~\ref{thm:main-almost-param} is optimal in general, see \cite{Rom19}, and the condition that $X$ be locally geodesic can be weakened, see Remark~\ref{rem:loc-geod-weakened} below. Theorem~\ref{thm:main-almost-param} also implies an analog in which the Jordan domain $\Omega$ is replaced by any open simply connected subset of $X$ with compact closure, see Corollary~\ref{cor:almost-param-Riem-map-version}.

We now describe conditions on $X$ under which a map $u$ as in the theorem can be upgraded to a homeomorphism, to a geometrically quasiconformal homeomorphism, or to a quasisymmetric homeomorphism. Given subsets $E,F,G\subset X$ we denote by $\MOD(E,F;G)$ the modulus of the family of curves joining $E$ and $F$ in $G$. Let $u$ be a map as in Theorem~\ref{thm:main-almost-param}. If for every $x\in X$ and every $R>0$ with $X\setminus B(x,R)\not=\emptyset$ we have 
\begin{equation}\label{eq:points-zero-modulus}
 \lim_{r\to 0} \MOD(B(x,r), X\setminus B(x,R); \bar{B}(x,R)) = 0
\end{equation}
then $u$ is a homeomorphism, see Proposition~\ref{prop:points-modulus-zero-implies-homeo}. If, moreover, there exists $\kappa>0$ such that every closed topological square $Q\subset X$ with boundary edges $\zeta_1,\zeta_2,\zeta_3,\zeta_4$ in cyclic order satisfies
\begin{equation}\label{eq:top-squares-moduli-opposite-edges}
 \MOD(\zeta_1, \zeta_3; Q)\cdot\MOD(\zeta_2,\zeta_4; Q)\leq \kappa
\end{equation}
then $u$ is geometrically quasiconformal, as follows from the arguments in \cite{Raj14}, see Proposition~\ref{prop:reciprocal-implies-second-modulus-inequality} below. Theorem~\ref{thm:main-almost-param} thus yields the following result, which is variant for locally geodesic metric surfaces of Rajala's uniformization theorem \cite[Theorem 12.1]{Raj14}.

\bc\label{cor:geom-qc-param-reciprocal}
 Let $X$ be a locally geodesic metric space homeomorphic to $\R^2$ and of locally finite Hausdorff $2$--measure. If $X$ satisfies conditions \eqref{eq:points-zero-modulus} and \eqref{eq:top-squares-moduli-opposite-edges} then for every Jordan domain $\Omega\subset X$ there exists a homeomorphism $u\colon\overline{D}\to \overline{\Omega}$ which is geometrically quasiconformal.
\ec

Rajala's reciprocality condition \cite{Raj14} mentioned above consists of  \eqref{eq:points-zero-modulus} and \eqref{eq:top-squares-moduli-opposite-edges} together with a lower bound on the product in \eqref{eq:top-squares-moduli-opposite-edges}. It has recently been shown in \cite{RR19} that the lower bound is always satisfied, see also \cite{EBPC21}. Notice that reciprocality is also a necessary condition for the existence of geometrically quasiconformal homeomorphisms. This is because $\R^2$ satisfies \eqref{eq:points-zero-modulus} and \eqref{eq:top-squares-moduli-opposite-edges} and they are preserved under geometrically quasiconformal homeomorphisms. The optimal geometric quasiconformality constants were obtained in \cite{Rom19}.
It was shown in \cite[Theorem 1.6]{Raj14} that if there exists $C>0$ such that 
\begin{equation}\label{eq:measure-balls-quadratic}
 \hm^2(B(x,r))\leq Cr^2
\end{equation}
for every $x\in X$ and $r>0$ then $X$ is reciprocal. It follows from \cite[Theorem 2.5]{LW20-param} that if $u$ is a map as in Theorem~\ref{thm:main-almost-param} and if $X$ satisfies \eqref{eq:measure-balls-quadratic} then $u$ is a quasisymmetric homeomorphism if and only if $\overline{\Omega}$ is linearly locally connected. In particular, Theorem~\ref{thm:main-almost-param} recovers the Bonk-Kleiner quasisymmetric uniformization theorem \cite{BK02} for metric discs (see Corollary~\ref{cor:qs-param} below) and, by a quasisymmetric gluing argument as in \cite{LW20-param}, also for metric spheres.

\subsection{Methods of proof and other results}

We now describe the ingredients in the proof of Theorem~\ref{thm:main-almost-param}. As in the classical existence proof of conformal parametrizations of smooth surfaces, we will obtain a quasiconformal almost parametrization $u$ as an energy minimizing disc in $X$ spanning $\partial \Omega$. This is similar to the approach taken in \cite{LW20-param}. The proofs in \cite{LW20-param} heavily use regularity properties and the intrinsic structure of energy minimizers in spaces with a quadratic isoperimetric inequality established in \cite{LW15-Plateau} and \cite{LW16-intrinsic}. Such results are not available in our setting.

Let $X$ be a complete metric space and let $N^{1,2}(D,X)$ be the space of Newton-Sobolev maps from $D$ to $X$ in the sense of \cite{HKST15}. For a map $u\in N^{1,2}(D, X)$ we denote by $\trace(u)$ the trace of $u$ and by $E_+^2(u)$ its (Reshetnyak) energy. If $\Gamma\subset X$ is a Jordan curve we denote by $\Lambda(\Gamma, X)$ the possibly empty family of maps $u\in N^{1,2}(D,X)$ whose trace has a continuous representative which is a weakly monotone parametrization of $\Gamma$. See Section~\ref{sec:prelims} for the definitions of these concepts.

\bt\label{thm:cont-energy-min-intro}
 Let $X$ be a locally geodesic metric space homeomorphic to $\overline{D}$ and let $\Gamma\subset X$ be a Jordan curve. If $u\in\Lambda(\Gamma, X)$ minimizes the Reshetnyak energy $E_+^2$ among all maps in $\Lambda(\Gamma, X)$ then $u$ has a representative which is continuous and extends continuously to the boundary.
\et

Notice that the regularity results for energy minimizers proved in \cite{LW15-Plateau} cannot be applied here since metric spaces as in Theorem~\ref{thm:cont-energy-min-intro} need not admit a quadratic isoperimetric inequality for curves. In general, the family $\Lambda(\Gamma, X)$ in Theorem~\ref{thm:cont-energy-min-intro} can be empty. However, we can prove the following theorem.

\bt\label{thm:non-trivial-Sobolev-Jordan-intro}
 Let $X$ be a locally geodesic metric space homeomorphic to $\R^2$, $\overline{D}$, or $S^2$. If $X$ has locally finite Hausdorff $2$--measure then $\Lambda(\Gamma, X)\not=\emptyset$ for every rectifiable Jordan curve $\Gamma\subset X$.
\et

 Theorem~\ref{thm:main-almost-param} now easily follows from Theorems~\ref{thm:cont-energy-min-intro} and \ref{thm:non-trivial-Sobolev-Jordan-intro} together with the results on the existence and structure of area and energy minimizers established in \cite{LW15-Plateau} and \cite{LW20-param}. Indeed, one easily reduces the theorem to the special case that $X$ is geodesic, homeomorphic to $\overline{D}$, and $\overline{\Omega}= X$. Since the boundary circle $\partial X$ of $X$ has finite length,  Theorem~\ref{thm:non-trivial-Sobolev-Jordan-intro} shows that the family $\Lambda(\partial X, X)$ is not empty. By \cite{LW15-Plateau} there exists an energy minimizer $u$ in $\Lambda(\partial X, X)$. By Theorem~\ref{thm:cont-energy-min-intro}, any such $u$ (has a representative which) is continuous up to the boundary and it thus follows from \cite[Theorem 1.2]{LW20-param} that $u$ is monotone. Finally, by \cite{LW15-Plateau} energy minimizers are infinitesimally quasiconformal in the sense that 
 \begin{equation}\label{eq:analytic-qc-ineq}
  (g_u(z))^2\leq K \cdot\jac(\apmd u_z)
 \end{equation}
 with $K=\frac{4}{\pi}$ for almost every $z\in D$, where $g_u$ denotes the minimal weak upper gradient of $u$ and $\jac(\apmd u_z)$ is the jacobian of $u$, see Section~\ref{sec:prelims}. This implies that $u$ satisfies \eqref{eq:analytic-qc-mod}, see Section~\ref{sec:properties-qc-almost-parametrizations}, and thus the outline of the proof of Theorem~\ref{thm:main-almost-param} is complete. Notice that a homeomorphism $u\colon \overline{D}\to X$ is quasiconformal according to the analytic definition if $u$ belongs to $N^{1,2}(D,X)$ and satisfies \eqref{eq:analytic-qc-ineq}, see Section~\ref{sec:properties-qc-almost-parametrizations}.

We mention the following consequence of Theorem~\ref{thm:non-trivial-Sobolev-Jordan-intro}.

\bc\label{cor:intro-rectifability}
 Let $X$ be a locally geodesic metric space homeomorphic to $\R^2$. If $X$ has locally finite Hausdorff $2$--measure then $X$ contains a $2$--rectifiable subset of positive Hausdorff measure.
\ec

An example showing that $X$ need not be countably $2$--rectifiable is given in \cite[Theorem A.1]{SWS10}. Corollary~\ref{cor:intro-rectifability} also holds for compact metric spaces of any topological dimension $n$ and finite Hausdorff $n$--measure with positive lower density almost everywhere. This was proved in \cite{DLD20} using a deep result of Bate \cite{Bat20} about purely unrectifiable metric spaces.

The paper is structured as follows. In Section~\ref{sec:prelims} we collect the necessary definitions and some results on Newton-Sobolev maps that will be used later. In  Section~\ref{sec:properties-qc-almost-parametrizations} we show how quasiconformal almost parametrizations can be upgraded under additional conditions on the underlying space. We prove Theorem~\ref{thm:cont-energy-min-intro} in Section~\ref{sec:cont-energy-min}; Theorem~\ref{thm:non-trivial-Sobolev-Jordan-intro} and Corollary~\ref{cor:intro-rectifability} are established in Section~\ref{sec:non-trivial-Sobolev-maps}. In the final Section~\ref{sec:proofs-almost-param} we discuss the proofs of the parametrizations results as well as some consequences and generalizations.

After having completed our manuscript we were informed that Dimitrios Ntalampekos and Matthew Romney are in the process of finalizing a paper \cite{NR21} in which they also prove (a variant of) Theorem~\ref{thm:main-almost-param} with a different approach.

\bigskip

{\bf Acknowledgments:} The second author would like to thank Kai Rajala for many inspiring discussions on the uniformization problem.

\section{Preliminaries}\label{sec:prelims}

\subsection{Basic notation}
We denote the open and closed unit discs in the Euclidean plane $\R^2$ by $D$ and $\overline{D}$, respectively; that is, 
$$D\coloneqq \left\{z\in\R^2: |z|<1\right\},\quad \overline{D}\coloneqq \left\{z\in\R^2: |z|\leq 1\right\},$$ where $|v|$ denotes the Euclidean norm of the vector $v\in\R^2$.
Let $(X,d)$ be a metric space. The open and closed balls in $X$ centered at some point $x_0$ of radius $r>0$ are $$B(x_0, r)\coloneqq \{x\in X: d(x,x_0)<r\},\quad \bar{B}(x_0,r)\coloneqq \{x\in X: d(x,x_0)\leq r\}.$$

Let $c\colon I\to X$ be a curve defined on some interval $I\subset\R$. The length of $c$ is denoted by $\length(c)$. If $c$ is absolutely continuous then $c$ has a metric derivative almost everywhere, thus the limit $$|c'(t)|:= \lim_{s\to t}\frac{d(c(s),c(t))}{|s-t|}$$ exists for almost every $t\in I$, and moreover $\length(c) = \int_I |c'(t)|\,dt$, see \cite{Kir94}.
A curve $c\colon[a,b]\to X$ is called geodesic if $\length(c) = d(c(a),c(b))$. The metric space $X$ is called geodesic if any pair of points in $X$ can be joined by a geodesic. It is called locally geodesic if every point $x\in X$ has a neighborhood $U$ such that any two points in $U$ can be joined by a geodesic in $X$.

Given $m\geq 0$, the $m$-dimensional Hausdorff measure on $X$ is denoted by $\hm^m$. The normalizing constant is chosen so that $\hm^n$ agrees with the Lebesgue measure on Euclidean $\R^n$. We write $|A|$ for the Lebesgue measure of a subset $A\subset\R^n$.

\subsection{Conformal modulus}

Let $X$ be a metric space and $\Gamma$ a family of curves in $X$. A Borel function $\rho\colon X\to [0,\infty]$ is said to be admissible for $\Gamma$ if $\int_\gamma \rho\geq 1$ for every locally rectifiable curve $\gamma\in\Gamma$. See \cite{HKST15} for the definition of the path integral $\int_\gamma\rho$.
The modulus of $\Gamma$ is defined by $$\MOD(\Gamma):= \inf_\rho\int_X\rho^2\,d\hm^2,$$ where the infimum is taken over all admissible functions for $\Gamma$. We emphasize that throughout this paper, the reference measure on $X$ will always be the $2$--dimensional Hausdorff measure. By definition, $\MOD(\Gamma)=\infty$ if $\Gamma$ contains a constant curve. A property is said to hold for almost every curve in $\Gamma$ if it holds for every curve in $\Gamma_0$ for some family $\Gamma_0\subset \Gamma$ with $\MOD(\Gamma\setminus \Gamma_0)=0$. In the definition of $\MOD(\Gamma)$, the infimum can equivalently be taken over all weakly admissible functions, that is, Borel functions $\rho\colon X\to [0,\infty]$ such that $\int_\gamma \rho\geq 1$ for almost every locally rectifiable curve $\gamma\in\Gamma$.

\subsection{Metric space valued Sobolev maps}

We recall some definitions from the theory of metric space valued Sobolev maps based on upper gradients \cite{Shan00}, \cite{HKST01}, \cite{HKST15} as well as two results concerning the existence and structure of energy minimizing discs established in \cite{LW15-Plateau}, \cite{LW20-param}. Note that the results in \cite{LW15-Plateau} are stated using a different but equivalent definition of Sobolev mappings.

Let $(X,d)$ be a complete metric space and $U\subset \R^2$ a bounded domain.  A Borel function $g\colon U\to [0,\infty]$ is said to be an upper gradient of a map $u\colon U\to X$ if
\begin{equation}\label{eq:upper-grad}
 d(u(\gamma(a)), u(\gamma(b)))\leq \int_\gamma g
\end{equation}
 for every rectifiable curve $\gamma\colon[a,b]\to U$. If \eqref{eq:upper-grad} only holds for almost every curve $\gamma$ then $g$ is called a weak upper gradient of $u$. A weak upper gradient $g$ of $u$ is called minimal weak upper gradient of $u$ if $g\in L^2(U)$ and if for every weak upper gradient $h$ of $u$ in $L^2(U)$ we have $g\leq h$ almost everywhere on $U$.

Denote by $L^2(U, X)$ the collection of measurable and essentially separably valued maps $u\colon U\to X$ such that the function $u_x(z):= d(u(z), x)$ belongs to $L^2(U)$ for some and thus any $x\in X$. A map $u\in L^2(U, X)$ belongs to the Newton-Sobolev space $N^{1,2}(U, X)$ if $u$ has a weak upper gradient in $L^2(U)$. Every such map $u$ has a minimal weak upper gradient $g_u$, unique up to sets of measure zero, see \cite[Theorem 6.3.20]{HKST15}.
The {\it Reshetnyak energy} of a map $u\in N^{1,2}(U, X)$ is defined by $$E_+^2(u):= \|g_u\|_{L^2(U)}^2.$$

If $u\in N^{1,2}(U, X)$ then for almost every $z\in U$ there exists a unique semi-norm on $\R^2$, denoted by $\apmd u_z$ and called the approximate metric derivative of $u$, such that 
$$\ap \lim _{y\to z}  \frac {d(u(y),u(z))- \apmd u_z(y-z)} {|y-z|} =0,$$ see \cite{Kar07} and \cite[Proposition 4.3]{LW15-Plateau}. For the definition of the approximate limit $\ap\lim$ see \cite{EG92}. The following notion of parametrized area was introduced in \cite{LW15-Plateau}.

\bd
 The parametrized (Hausdorff) area of a map $u\in N^{1,2}(U, X)$ is defined by $$\Area(u)= \int_U \jac(\apmd u_z)\,dz,$$ where the Jacobian $\jac(s)$ of a semi-norm $s$ on $\R^2$ is the Hausdorff $2$-measure on $(\R^2, s)$ of the unit square if $s$ is a norm and zero otherwise.
 \ed
 
The area of the restriction of $u$ to a measurable set $B\subset U$ is defined analogously. It is not difficult to show that $\jac(\apmd u_z)\leq (g_u(z))^2$ for almost every $z\in U$, see \cite[Lemma 7.2]{LW15-Plateau}. If $u$ is a homeomorphism onto its image then the jacobian $\jac(\apmd u_z)$ agrees with the Radon-Nikodym derivative of the measure $u^*\hm^2(B):= \hm^2(u(B))$ with resepct to the Lebesgue measure at almost every point $z\in U$.

\bd\label{def:inf-qc}
 A map $u\in N^{1,2}(U, X)$ is called infinitesimally $K$-quasiconformal if
\begin{equation}\label{eq:inf-qc}
(g_u(z))^2 \leq K\cdot \jac(\apmd u_z)
\end{equation}
for almost every $z\in U$. 
\ed

If $u\in N^{1,2}(D, X)$ then for almost every $v\in S^1$ the curve $t\mapsto u(tv)$ with $t\in[1/2, 1)$ is absolutely continuous. The trace of $u$ is defined by $\trace(u)(v):= \lim_{t\nearrow 1}u(tv)$ for almost every $v\in S^1$. It follows from \cite{KS93} that $\trace(u)\in L^2(S^1, X)$. If $u$ is the restriction to $D$ of a continuous map $\hat{u}$ on $\overline{D}$ then $\trace(u)=\hat{u}|_{S^1}$. 
Given a Jordan curve $\Gamma\subset X$ we denote by $\Lambda(\Gamma, X)$ the possibly empty family of maps $v\in N^{1,2}(D, X)$ whose trace has a continuous representative which weakly monotonically parametrizes $\Gamma$. Recall that a continuous map $c\colon S^1\to \Gamma$ is called a weakly monotone parametrization of $\Gamma$ if $c$ is the uniform limit of homeomorphisms $c_i\colon S^1\to \Gamma$.

\bt\label{thm:existence-energy-min}
 Let $X$ be a proper metric space and $\Gamma\subset X$ be a Jordan curve. If $\Lambda(\Gamma, X)$ is not empty then there exists $u\in \Lambda(\Gamma, X)$ satisfying $$E_+^2(u) = \inf\left\{E_+^2(v): v\in\Lambda(\Gamma, X)\right\},$$ and any such $u$ is infinitesimally $K$--quasiconformal with $K=\frac{4}{\pi}$.
\et

\begin{proof}
The existence of an energy minimizer in $\Lambda(\Gamma, X)$ follows from \cite[Theorem 7.6]{LW15-Plateau}. Energy minimizers are infinitesimally $K$--quasiconformal with $K=\frac{4}{\pi}$ by \cite{LW16-energy-area}, see also \cite[Lemma 6.5]{LW15-Plateau}.
\end{proof}

We will also need the following theorem proved in \cite{LW20-param}.

\bt\label{thm:energy-min-monotone}
 Let $X$ be a geodesic metric space homeomorphic to $\overline{D}$, and let $u\colon\overline{D}\to X$ be a continuous map. If $u$ belongs to $\Lambda(\partial X, X)$ and minimizes the Reshetnyak energy $E_+^2$ among all maps in $\Lambda(\partial X, X)$ then $u$ is monotone.
\et

By definition, the boundary circle $\partial X$ of $X$ is the image of $S^1$ under a homeomorphism from $\overline{D}$ to $X$. Recall that a continuous map $u\colon \overline{D}\to X$ is monotone if $u^{-1}(x)$ is connected for every $x\in X$. If $X$ is homeomorphic to $\overline{D}$ then $u$ is monotone if and and only if $u$ is the uniform limit of homeomorphisms $u_n\colon \overline{D}\to X$, see \cite{You48}.

\section{Upgrading a quasiconformal almost parametrization}\label{sec:properties-qc-almost-parametrizations}

The aim of this short section is to summarize some results which show that maps as in Theorem~\ref{thm:main-almost-param} can be upgraded under certain additional conditions on the underlying space.

We first recall the connection with infinitesimally quasiconformal maps. 
Let $X$ be a complete metric space and $u\colon \overline{D}\to X$ continuous and monotone. If $u\in N^{1,2}(D, X)$ and $u$ is infinitesimally $K$--quasiconformal then 
\begin{equation}\label{eq:inf-qc-props}
    \MOD(\Gamma)\leq K\cdot\MOD(u\circ\Gamma)
\end{equation}
for every family $\Gamma$ of curves in $\overline{D}$. Conversely, if $u$ is a homeomorphism onto its image and satisfies \eqref{eq:inf-qc-props} then $u$ belongs to $N^{1,2}(D, X)$ and is infinitesimally $K$--quasiconformal. See \cite[Proposition 3.5]{LW20-param} and \cite[Theorem 1.1]{Wil12} for a proof.

\bp\label{prop:points-modulus-zero-implies-homeo}
 Let $X$ be a complete metric space satisfying \eqref{eq:points-zero-modulus}. Let $u\colon \overline{D}\to X$ be continuous, monotone, and non-constant. If $u$ satisfies \eqref{eq:inf-qc-props} then $u$ is a homeomorphism onto its image.
\ep

\begin{proof}
 This follows exactly as in the proof of \cite[Theorem 3.6]{LW20-param}.
\end{proof}

In the setting of Theorem~\ref{thm:main-almost-param} there need not exist a homeomorphism satisfying \eqref{eq:inf-qc-props}. The following example illustrating this appears in \cite[Example 11.3]{LW16-intrinsic}, see \cite{Raj14} for other examples.

\begin{example}\label{ex:collapsed-disc}{\rm
 Let $T=\{z\in D: |z|\leq 1/2\}$ and let $X= \overline{D}/T$ be the quotient metric space equipped with the quotient metric. Then $X$ is geodesic, homeomorphic to $\overline{D}$, and has finite Hausdorff $2$--measure. We claim that there does not exist a homeomorphism $u\colon \overline{D}\to X$ satisfying \eqref{eq:inf-qc-props}. Suppose to the contrary that such $u$ exists. Then $u$ is analytically quasiconformal by the discussion above, thus $u$ is in $N^{1,2}(D, X)$ and is infinitesimally quasiconformal. Let $\pi\colon \overline{D}\to X$ be the quotient map and set $p:=\pi(T)$. After possibly precomposing $u$ with a biLipschitz homeomorphism of $\overline{D}$ we may assume that $u(0)=p$. Consider the homeomorphism $v\colon D\setminus\{0\} \to D\setminus T$ satisfying $\pi(v(z)) = u(z)$ for all $z\in D\setminus\{0\}$. Since $\pi$ is a local isometry on $D\setminus T$ it follows that $v$ is (analytically) quasiconformal, which is impossible since the punctured disc is not quasiconformally equivalent to the annulus, see \cite[Theorem 39.1]{Vai71}. This contradiction finishes the proof of the claim.}
\end{example}

The next proposition follows from the arguments in \cite[Section 11]{Raj14}.

\bp\label{prop:reciprocal-implies-second-modulus-inequality}
 Let $X$ be a metric space homeomorphic to $\overline{D}$. Suppose $u\colon\overline{D}\to X$ is a homeomorphism satisfying
\eqref{eq:inf-qc-props}. Then $u$ is geometrically quasiconformal if and only if $X$ satisfies \eqref{eq:top-squares-moduli-opposite-edges} for some $\kappa$.
\ep

\begin{proof}
 Notice that $\overline{D}$ satisfies \eqref{eq:top-squares-moduli-opposite-edges}. Therefore, if $u$ is geometrically quasiconformal then also $X$ satisfies  \eqref{eq:top-squares-moduli-opposite-edges}. Suppose now that $X$ satisfies \eqref{eq:top-squares-moduli-opposite-edges} for some $\kappa$.
 By the discussion at the beginning of this section, the map $u$ thus belongs to $N^{1,2}(D,X)$ and is infinitesimally $K$-quasiconformal. Identifying $\overline{D}$ with the unit square $R:=[0,1]^2$ via a biLipschitz homeomorphism, we may view $u$ as an element of $N^{1,2}(R,X)$. There exists a Borel set $A\subset R$ of full measure such that $u|_A$ has Lusin's property (N), see e.g.~\cite[Proposition 3.2]{LW15-Plateau}. Let $g_u$ be the weak minimal upper gradient of $u$. We may assume that $g_u=\infty$ on $R\setminus A$. The Borel function $h\colon X\to[0,\infty]$ defined by $h:= \frac{1}{g_u\circ u^{-1}}$ is $L^2$--integrable since $$\int_X h^2\,d\hm^2 = \int_{u(A)} h^2\,d\hm^2 = \int_A \frac{1}{g_u^2(z)}\jac(\apmd u_z)\,dz\leq |R|<\infty$$ by the area formula, see \cite{Kir94} and \cite{Kar07}. Arguing exactly as in the proof of \cite[Proposition 11.1]{Raj14} one shows that there exists $C$ only depending on $K$ and $\kappa$ such that $C\cdot h$ is a weak upper gradient for $u^{-1}$. The proof of this relies on a lower bound of the form
  \begin{equation}\label{eq:lower-bound-int-h-Q}
  \int_{Q(i,j,k)} h^2\,d\hm^2 \geq K^{-1}2^{-2k},
 \end{equation}
 where $Q(i,j,k) = u([2^{-k}i, 2^{-k}(i+1)]\times[2^{-k}j,2^{-k}(j+1)])$
 as well as on an upper bound of the form
 \begin{equation}\label{eq:upper-bound-mod-Gamma-ell}
  \MOD(\Gamma_\ell(i,j,k))\leq 3\kappa
 \end{equation}
 for suitable path families $\Gamma_\ell(i,j,k)$ defined in the proof of \cite[Proposition 11.1]{Raj14}. In our case, \eqref{eq:lower-bound-int-h-Q} follows from the area formula and \eqref{eq:upper-bound-mod-Gamma-ell} follows from \eqref{eq:inf-qc-props} and \eqref{eq:top-squares-moduli-opposite-edges}. 
 Now let $\Gamma$ be a family of curves in $R$ and let $\varrho$ be an admissible function for $\Gamma$. Since for almost every curve $\beta = u\circ\gamma\in u\circ\Gamma$ we have $$1\leq \int_\gamma \varrho\leq C\cdot\int_{\beta}\varrho\circ u^{-1} h$$ it follows that $C\varrho\circ u^{-1} h$ is weakly admissible for $u\circ\Gamma$. The area formula yields
 \begin{equation*}
      \MOD(u\circ\Gamma)\leq C^2\cdot\int_Xh^2\varrho^2\circ u^{-1}\,d\hm^2
      \leq C^2\int_R\varrho^2(z)\,dz,
 \end{equation*}
 and taking the infimum over $\varrho$ we conclude that $\MOD(u\circ\Gamma)\leq C^2\cdot\MOD(\Gamma)$.
\end{proof}

We finally describe conditions that imply that $u$ is quasisymmetric. Recall that a homeomorphism $u\colon M\to X$ between metric spaces is said to be quasisymmetric if there exists a homeomorphism $\eta\colon[0,\infty)\to[0,\infty)$ such that $$d(u(z),u(a))\leq \eta(t)\cdot d(u(z),u(b))$$ for all $z,a,b\in M$ with $d(z,a)\leq t\cdot d(z,b)$. 

\bp\label{prop:upgrade-qs}
 Let $X$ be a metric space homeomorphic to $\overline{D}$ and let $u\colon \overline{D}\to X$ be a homeomorphism satisfying \eqref{eq:inf-qc-props}. If there exists $L>0$ such that $$\hm^2(B(x,r))\leq L\cdot r^2$$ for all $x\in X$ and $r>0$ then $u$ is quasisymmetric if and only if $X$ is linearly locally connected.
\ep

We refer for example to the Appendix of \cite{LW20-param} for a proof of the proposition.
Here, a metric space $X$ is called linearly locally connected if there exists $\lambda\geq 1$ such that for every $x\in X$ and for all $r>0$, every pair of points in $B(x,r)$ can be joined by a continuum in $B(x,\lambda r)$ and every pair of points in $X\setminus B(x,r)$ can be joined by a continuum in $X\setminus B(x,r/\lambda)$.

\section{Continuity of energy minimizers in locally geodesics metric discs}\label{sec:cont-energy-min}

In this section we prove Theorem~\ref{thm:cont-energy-min-intro}. For an arbitrary map $v\colon D\to X$ to a metric space $X$ and for $z\in D$ and $\delta>0$ set $$\osc(v,z,\delta):= \inf\{\diam(v(A)): \text{$A\subset D\cap B(z,\delta)$ subset of full measure}\},$$ called the essential oscillation of $v$ in the $\delta$--ball around $z$.

\bp\label{prop:essential-oscillation-energy-min}
 Let $X$ be a locally geodesic metric space homeomorphic to $\overline{D}$ and let $\Gamma\subset X$ be a Jordan curve. Suppose $u\in\Lambda(\Gamma,X)$ minimizes the Reshetnyak energy among all maps  in $\Lambda(\Gamma, X)$. Then for every $\varepsilon>0$ there exists $\delta>0$ such that $\osc(u,z,\delta)<\varepsilon$ for every $z\in D$.
\ep

The theorem easily follows from this proposition.

\begin{proof}[Proof of Theorem~\ref{thm:cont-energy-min-intro}]
 Let $A=\{z_n: n\in\N\}\subset D$ be a countable dense set. For each $k\in\N$ apply the proposition with $\varepsilon = \frac{1}{k}$ to obtain $\delta_k>0$ and negligible subsets $N_{k,n}\subset D$ such that $$\diam(u(D\cap B(z_n,\delta_k)\setminus N_{k,n}))<\frac{1}{k}$$ for all $n\in\N$. Then the set $N:= \cup_{k,n\in\N} N_{k,n}$ is negligible. Let $\varepsilon>0$ and choose $k\in\N$ such that $\frac{1}{k}<\varepsilon$. If $z,z'\in D\setminus N$ satisfy $|z-z'|<\delta_k$ then there exists $n$ such that $z,z'\in B(z_n,\delta_k)\setminus N_{k,n}$ and hence $$d(u(z), u(z'))\leq \diam(u(D\cap B(z_n,\delta_k)\setminus N_{k,n}))<\frac{1}{k} <\varepsilon.$$ This shows that $u|_{D\setminus N}$ is uniformly continuous and hence has a unique continuous extension $\bar{u}$ to $\overline{D}$.
\end{proof}

We need the following lemma in the proof of Proposition~\ref{prop:essential-oscillation-energy-min}.

\bl\label{lem:biLip-Jordan}
 Let $X$ be a locally geodesic metric space homeomorphic to $\overline{D}$ and let $\varepsilon>0$. Then there is $\varepsilon'>0$ such that for every $x\in X$ there exists a biLipschitz curve $T\subset X$ with the following property. Either $T$ is the boundary of a Jordan domain $\Omega$ containing  $\bar{B}(x,\varepsilon')$ and with $\diam\Omega\leq \varepsilon$; or $T$ is a Jordan arc intersecting $\partial X$ exactly at its endpoints, and a component $\Omega$ of $X\setminus T$ contains $\bar{B}(x,\varepsilon')$ and satisfies $\diam \Omega\leq \varepsilon$.
\el

\begin{proof}
 Let $0<\varepsilon<\diam X$ and let $\varrho\colon\overline{D}\to X$ be a homeomorphism. Choose $\varepsilon', \delta>0$ such that $$B(\varrho(z), 2\varepsilon')\subset \varrho(\overline{D}\cap B(z,\delta)) \subset B(\varrho(z), \varepsilon/3)$$ for every $z\in\overline{D}$. Let $x\in X$, set $z=\varrho^{-1}(x)$ and $S=\{w\in\overline{D}: |z-w|=\delta\}$. We can approximate the curve $\varrho(S)$ by a biLipschitz curve $T\subset X$ with the desired properties as follows. 
 
 We distinguish two cases and first assume that $S$ does not intersect the boundary of $\overline{D}$ and hence is a circle. Let $\alpha\colon S^1\to S$ be a parametrization and let $r>0$ be sufficiently small, to be determined later. Let $\{t_0, t_1,\dots, t_n\}$ be a fine partition of $S^1$ and let $\gamma\colon S^1\to X$ be a curve such that $\gamma(t_k) = \varrho(\alpha(t_k))$ and the restriction of $\gamma$ to the short segment $I_k=\overline{t_k t_{k+1}}\subset S^1$ is a geodesic for every $k$. If the partition is chosen sufficiently fine, we have $\diam(\varrho^{-1}(\gamma(I_k))) <\frac{r}{8}$ for all $k$. In particular, if $s, t\in S^1$ are such that $\gamma(s) = \gamma(t)$ then the shorter segment $\overline{st}\subset S^1$ is such that $\varrho^{-1}(\gamma(\overline{st}))$ is contained in a ball of radius $r/2$ centered on $S$ and thus homotopic relative to endpoints to the constant curve inside this ball. It is then not difficult to see that, after deleting a finite number of subcurves from $\gamma$, we obtain a piecewise geodesic Jordan curve $T\subset X$ such that $\varrho^{-1}(T)$ is homotopic to $S$ in the $r$--neighborhood of $S$. Moreover, by applying the claim in the proof of \cite[Lemma 4.2]{LW20-param} we may further assume that $T$ is a biLipschitz Jordan curve. The Jordan domain $\Omega'$ enclosed by $\varrho^{-1}(T)$ satisfies $B(z,\delta - r)\subset \Omega'$. Thus, if $r>0$ was chosen small enough then the Jordan domain $\Omega = \varrho(\Omega')$ enclosed by $T$ satisfies $\bar{B}(x,\varepsilon')\subset\Omega \subset B(x, \varepsilon/2)$, as desired.
 
 The case that $S$ intersects the boundary of $\overline{D}$ is analogous and is left to the reader.
\end{proof}

We now prove Proposition~\ref{prop:essential-oscillation-energy-min}. For this, let $u\in\Lambda(\Gamma, X)$ be an energy minimizer and denote by $\alpha\colon S^1\to X$ the continuous representative of the trace of $u$. Let $\varepsilon>0$ and let $\varepsilon'>0$ be as in Lemma~\ref{lem:biLip-Jordan}. Choose $\delta\in(0,1/4)$ so small that $2\pi^2 E_+^2(u)<(\varepsilon')^2|\log\delta|$ and such that $\alpha$ maps segments of diameter $2\sqrt{\delta}$ to sets of diameter at most $\varepsilon'$. 

Fix $z\in D$ and for $r>0$ let $\gamma_r$ be a constant speed parametrization of the curve $\{p\in \overline{D}:|p-z|=r\}$. By the Courant-Lebesgue Lemma (see e.g. \cite[Lemma~7.3]{LW15-Plateau}) there exists a set $A\subset(\delta,\sqrt{\delta})$ of strictly positive measure such that $u\circ\gamma_r$ has an absolutely continuous representative, denoted again by $u\circ\gamma_r$, of length $$\length(u\circ\gamma_r)\leq \pi\left(\frac{2E_+^2(u)}{|\log\delta|}\right)^{\frac{1}{2}}< \varepsilon'$$ for every $r\in A$.
For almost every $r\in A$ for which $\gamma_r$ intersects $S^1$, the endpoints of the absolutely continuous curve $u\circ\gamma_r$ coincide with $\alpha(a_r)$ and $\alpha(b_r)$, where $a_r,b_r$ are the endpoints of $\gamma_r$. We furthermore have $$\trace(u|_{D\cap B(z,r)})\circ\gamma_r = u\circ\gamma_r$$ for almost every $r\in A$. 
Fix $r\in A$ such that all of the above hold and set $W:= D\cap B(z,r)$. Since $\alpha$ maps segments of diameter $2\sqrt{\delta}$ to sets of diameter at most $\varepsilon'$ it follows that the image of (the continuous representative of) the trace of $u|_W$ is contained in $\bar{B}(x,\varepsilon')$ for some $x\in X$. By Lemma~\ref{lem:biLip-Jordan}, there exists a biLipschitz curve $T\subset X$ and a set $\Omega$ with $\diam(\Omega)\leq\varepsilon$ and $\bar{B}(x,\varepsilon')\subset\Omega$ such that $\Omega$ is either a Jordan domain and $T=\partial\Omega$ or $T$ is a Jordan arc intersecting $\partial X$ exactly at its endpoints and $\Omega$ is a component of $X\setminus T$. We now claim:

\bl\label{lem:negligible-set-mapped-outside}
 The set $N= \{w\in W: u(w)\in X\setminus\overline{\Omega}\}$ is negligible.
\el

With this lemma at hand we can easily finish the proof of the proposition. Indeed, we have $u(W\setminus N)\subset \overline{\Omega}$ and therefore $$\diam(u(W\setminus N))\leq \diam(\overline{\Omega})\leq \varepsilon.$$ Hence, since $|N|=0$, we obtain $$\osc(u,z,\delta)\leq \osc(u,z,r)\leq \varepsilon.$$ Since $\delta>0$ was independent of $z$ this proves the proposition.

We are left to prove the lemma above.

\begin{proof}[Proof of Lemma~\ref{lem:negligible-set-mapped-outside}]
 Since $u$ is an energy minimizer it follows that $u$ is infinitesimally quasiconformal and, by \cite[Theorem~1.1]{LW16-energy-area}, minimizes the inscribed Riemannian area
    $$\Area_{\mu^i}(u):=\int_D\jac_{\mu^i}(\apmd u_z)dz$$
    among all maps in $\Lambda(\Gamma, X)$. Here, the $\mu^i$-jacobian $\jac_{\mu^i}(s)$ of a semi-norm $s$ on $\R^2$ is given by $\jac_{\mu^i}(s)=0$ if $s$ is degenerate and $\jac_{\mu^i}(s)=\frac{\pi}{|L|}$ if $s$ is a norm, where $|L|$ denotes the Lebesgue measure of John's ellipse of $\{v\in\R^2:s(v)\leq1\}$.
    
 In order to prove that $N$ is negligible we suppose to the contrary that $|N|>0$. We then claim that $\Area_{\mu^i}(u|_N)>0$. In order to prove the claim, we argue as in the proof of \cite[Proposition 7]{GHP19} and decompose $W$ into horizontal curves $\beta_t$. For almost every $t$ the composition $u\circ\beta_t$ has an absolutely continuous representative $c_t$ with speed
    $$|c_t'(s)|=\apmd u_{\beta_t(s)}(\beta_t'(s))$$
    for almost every $s$, see \cite[Lemma 4.9]{LW15-Plateau}. If $\hm^1(\beta_t\cap N)>0$ then the set $c_t^{-1}(X\setminus\overline{\Omega})$ is nonempty and open and thus an at most countable disjoint union of open intervals. Almost every point in such an interval $I$ is contained in $N$ and, as $c_t|_I$ has endpoints in $\overline{\Omega}$, it follows that $\length(c_t|_I)>0$. We conclude
    $$\int_I\apmd u_{\beta_t(s)}(\beta_t'(s))\,ds =\int_I|c_t'(s)|\,ds= \length(c_t|_I)>0$$
    and therefore $\apmd u_{\beta_t(s)}(\beta_t'(s))$ can not vanish for almost every $s\in I$. By Fubini we thus obtain that $N$ contains a set $A$ of strictly positive measure such that $\apmd u_w\not= 0$ for every $w\in A$. Since $u$ is infinitesimally quasiconformal it follows that $\jac_{\mu^i}(\apmd u_w)>0$ for almost every $w\in A$  and hence $\Area_{\mu^i}(u|_N)>0$. This proves the claim.
    
 Next, let $T$ and $\Omega$ be as above and notice that there exists a continuous retraction $\varrho_0\colon X\to\overline{\Omega}$ such that $\varrho_0(X\setminus\Omega)\subset T$. Since $T$ is a biLipschitz curve, it is locally Lipschitz $n$--connected for every $n$. Moreover, $X$ has Nagata dimension at most $2$ by \cite{JL20}. Hence, by the Lipschitz extension results in \cite{Hoh93} and \cite{LS05}, we can approximate $\varrho_0$ arbitrarily closely by a Lipschitz retraction $\varrho\colon X \to \overline{\Omega}$ satisfying $\varrho(X\setminus \Omega) \subset T$. Let $v\colon D\to X$ be the map which agrees with $u$ on $D\setminus W$ and with $\varrho\circ u$ on $W$. Since the trace of $u|_W$ has image in $\Omega$ it follows from the gluing theorem \cite[Theorem 1.12.3]{KS93} that $v\in N^{1,2}(D, X)$ and $\trace(v) = \trace(u)$; in particular, we have $v\in\Lambda(\Gamma, X)$. Notice that $v$ agrees with $u$ on $D\setminus N$ and that $\Area_{\mu^i}(v|_N)=0$ since $v(N)\subset T$. It thus follows that $$\Area_{\mu^i}(u)=  \Area_{\mu^i}(u|_{D\setminus N}) + \Area_{\mu^i}(u|_N)> \Area_{\mu^i}(u|_{D\setminus N}) = \Area_{\mu^i}(v),$$ which contradicts the area minimization property of $u$. This completes the proof of the lemma.
\end{proof}

\section{Non-trivial Sobolev maps in locally geodesic metric surfaces}\label{sec:non-trivial-Sobolev-maps}

In this section we establish Theorem~\ref{thm:non-trivial-Sobolev-Jordan-intro}. In its proof we will use the fact that every compact metric space isometrically embeds into an injective metric space which is again compact. Recall that a metric space $E$ is injective if for every metric space $Z$, any $1$--Lipschitz map $Y\to E$, defined on a subset $Y\subset Z$, extends to a $1$--Lipschitz map $Z\to E$. By \cite{Isb64}, for every metric space $X$ there exists an injective metric space $E(X)$ which contains $X$ and which is minimal in an appropriate sense among injective metric spaces containing $X$. Such a space $E(X)$ is called injective hull of $X$ and is unique up to isometry. Moreover, if $X$ is compact then so is $E(X)$. See \cite{Isb64} for the proof of these properties.

The following proposition is the main ingredient in the proof of Theorem~\ref{thm:non-trivial-Sobolev-Jordan-intro}.

\bp\label{prop:Lipschitz-approx-neighborhood}
 Suppose $X$ is a geodesic metric space homeomorphic to $\overline{D}$. If $\hm^2(X)<\infty$ and $\length(\partial X)<\infty$ then there exists $M>0$ with the following property. For every $\varepsilon>0$ there is a Lipschitz map $v\colon \overline{D}\to E(X)$ with $\Area(v)\leq M$ and such that  $v|_{S^1}$ parametrizes $\partial X$ and the image of $v$ is contained in the $\varepsilon$--neighborhood of $X$ in $E(X)$.
\ep

It is not difficult to prove Theorem~\ref{thm:non-trivial-Sobolev-Jordan-intro} using this proposition.

\begin{proof}[Proof of Theorem~\ref{thm:non-trivial-Sobolev-Jordan-intro}]
 We first assume that $X$ is geodesic, homeomorphic to $\overline{D}$, and $\Gamma = \partial X$. By Proposition~\ref{prop:Lipschitz-approx-neighborhood} there exist $M>0$ and a sequence $(v_n)$ of Lipschitz maps $v_n\colon\overline{D}\to E(X)$ with $\Area(v_n)\leq M$ and such that $v_n|_{S^1}$ parametrizes $\partial X$ and the image of $v_n$ is contained in the $1/n$--neighborhood of $X$ for each $n\in\N$. 
 
 By Morrey's $\varepsilon$--conformality Lemma \cite{FW20} there exist diffeomorphisms $\varrho_n$ of $\overline{D}$ such that $u_n:= v_n\circ\varrho_n$ satisfies $$E_+^2(u_n)\leq \frac{4}{\pi}\cdot \Area(u_n) + 1 \leq \frac{4M}{\pi} +1$$ for every $n$. Let $p_1,p_2,p_3\in S^1$ and $q_1,q_2,q_3\in \partial X$ be distinct points. After precomposing $u_n$ with a Moebius transformation (this leaves the energy invariant) we may assume that every $u_n$ satisfies the $3$--point condition $u_n(p_i) = q_i$ for $i=1,2,3$. Thus, the sequence $(\alpha_n)$ of curves $\alpha_n:=u_n|_{S^1}$ is equi-continuous by \cite[Proposition 7.4]{LW15-Plateau}. Therefore, after passing to a subsequence, we may assume by the Arzel\`a-Ascoli theorem that $(\alpha_n)$ converges uniformly to a curve $\alpha$. As the uniform limit of parametrizations of $\partial X$, the curve $\alpha$ is a weakly monotone parametrization of $\partial X$. Finally, after passing to a further subsequence, we may assume by the Rellich-Kondrachov theorem \cite[Theorem 1.13]{KS93} that $(u_n)$ converges in $L^2(D,E(X))$ to some $u\in N^{1,2}(D, E(X))$. Since the image of $u_n$ is contained in the $1/n$--neighborhood of $X$ for every $n$ it follows that the essential image of $u$ is contained in $X$, so we may view $u$ as an element of $N^{1,2}(D, X)$. Since the traces converge in $L^2(S^1, E(X))$ to $\trace(u)$ by \cite[Theorem 1.12.2]{KS93} it follows that $\trace(u) = \alpha$ and hence that $u\in \Lambda(\partial X, X)$. This shows that $\Lambda(\partial X, X)$ is not empty, which completes the proof of the special case.
 
 Now, let $X$ and $\Gamma$ be as in the statement of the theorem. Then $\Gamma$ encloses a Jordan domain $\Omega\subset X$. Denote by $d$ the metric on $X$ and consider the length metric $d_{\overline{\Omega}}$ on $\overline{\Omega}$. The identity map $\pi\colon (\overline{\Omega},d_{\overline{\Omega}}) \to (\overline{\Omega}, d)$ is a homeomorphism, is $1$--Lipschitz, and preserves the lengths of curves and the Hausdorff $2$--measures of Borel subsets, see \cite[Lemma 2.1]{LW20-param}. In particular, the metric space $Y = (\overline{\Omega}, d_{\overline{\Omega}})$ is geodesic, homeomorphic to $\overline{D}$, and $\length(\partial Y)$ and $\hm^2(Y)$ are finite. It thus follows from the first part of the proof that $\Lambda(\partial Y, Y)$ is not empty. Let $v\in\Lambda(\partial Y, Y)$. Then $u\coloneqq \pi\circ v$ is an element of $N^{1,2}(D, X)$ with image in the compact set $\overline{\Omega}$. Since $\trace(u) = \pi\circ \trace(v)$ and $\trace(v)$ has a continuous representative which is a weakly monotone parametrization of $\partial Y$, we see that $u\in\Lambda(\Gamma, X)$. This shows that $\Lambda(\Gamma, X)$ is not empty and completes the proof.
\end{proof}

The remainder of this section is dedicated to the proof of  Proposition~\ref{prop:Lipschitz-approx-neighborhood}. We will need two lemmas.

\bl\label{lem:approximation-simplicial-complex}
 There is a constant $C\geq 1$ with the following property. Let $X$ be a geodesic metric space homeomorphic to $\overline{D}$. Then for every $r>0$ there exist a finite metric simplicial complex $\Sigma$ and $C$--Lipschitz maps $\psi\colon X\to \Sigma$ and $\varphi\colon \Sigma\to E(X)$ subject to:
 \begin{enumerate}
     \item[(i)] $\Sigma$ has dimension $\leq 2$ and the metric on $\Sigma$ is geodesic and such that every simplex is a Euclidean simplex of side length $r$;
     \item[(ii)] the image of $\varphi$ is in the $Cr$--neighborhood of $X$ and $d(x,\varphi(\psi(x)))\leq Cr$ for all $x\in X$.
 \end{enumerate}
\el

The lemma is a consequence of \cite[Theorem 2]{JL20} and \cite[Theorem 3.6]{BW21}. For the convenience of the reader we sketch the argument.

\begin{proof}
 The space $X$ has Nagata dimension at most $2$ with some universal constant $c$ by \cite[Theorem 2]{JL20}. Thus, for a given $r>0$, there exists a finite cover $\{B_1,\dots, B_k\}$ of $X$ by sets $B_i\subset X$ of diameter at most $cr$ and such that every subset of $X$ of diameter at most $r$ intersects at most three of the $B_i$'s. Define $1$--Lipschitz functions $\tau_i\colon X\to \R$ by $\tau_i(x) = \max\{\frac{r}{2}-d(x,B_i), 0\}$. Then for every $x$ we have $\bar{\tau}(x):= \tau_1(x)+\dots+\tau_k(x)\geq \frac{r}{2}$  and $\tau_i(x)>0$ for at most three indices $i$. Therefore, the map $\psi(x) = \bar{\tau}(x)^{-1}(\tau_1(x),\dots,\tau_k(x))$ has image in the $2$--skeleton of the simplex $\Delta= \{(v_1,\dots, v_k)\in\R^k: v_i\geq 0, v_1+\dots+v_k = 1\}$. One calculates as in the proof of \cite[Theorem 5.2]{LS05} that $$|\psi(x) - \psi(y)|\leq 24r^{-1}d(x,y)$$ for all $x,y\in X$. Let $\Sigma$ be the smallest subcomplex of $\Delta$ containing $\psi(X)$ and define a map $\varphi\colon\Sigma\to E(X)$ as follows. For each vertex $e_i\in\Sigma^{(0)}$ let $\varphi(e_i)$ be a point in $B_i$. If $e_i, e_j$ are adjacent vertices in $\Sigma$ then $d(\varphi(e_i),\varphi(e_j))\leq (2c+1)r|e_i-e_j|$. Using the Lipschitz connectedness of $E(X)$ we can extend $\varphi|_{\Sigma^{(0)}}$ to the $1$--simplices and $2$--simplices of $\Sigma$ and obtain a map $\varphi$ which is $Cr$--Lipschitz on each simplex and satisfies $d(x,\varphi(\psi(x)))\leq Cr$ for some $C$ only depending on $c$. Let $d_\Sigma$ be the length metric on $\Sigma$ and scaled by the factor $r/\sqrt{2}$. Then $(\Sigma, d_\Sigma)$, $\psi$, and $\varphi$ satisfy the properties of the lemma; see \cite[Section 3]{BW21} for details.
\end{proof}

Given sets $Y$, $Z$, and a map $f\colon Z\to Y$ we set for each $y\in Y$ $$N(f, y):=\#\{z\in Z: f(z) = y\},$$ the multiplicity of $f$ at $y$. If $Z$ is an open subset of $\R^n$ and $Y=\R^n$ and $f$ is continuous then the multiplicity function $N(f, \cdot)$ is (Lebesgue) measurable, see \cite{RR55}. We need the following lemma. See \cite[Lemma 7.3]{Fed55}, \cite{Fed48}, and \cite[Lemma A.1]{WY18} for related results.

\bl\label{lem:approx-cont-maps-multiplicity}
 Let $r>0$ and let $\Sigma$ be a finite simplicial complex of dimension at most $2$, equipped with a metric such that each simplex is a Euclidean simplex of side length $r$.
 Let furthermore $\varrho\colon\overline{D}\to \Sigma$ be a continuous map such that $\varrho(S^1)$ is contained in the $1$--skeleton $\Sigma^{(1)}$ of $\Sigma$ and $$L:= 4 \cdot 3^{-\frac{1}{2}} r^{-2} \int_\Sigma N(\varrho, y)\,d\hm^2(y)<\infty.$$ Then there exist disjoint compact balls $B_1,\dots, B_m\subset D$ for some $0\leq m\leq L$, and a continuous map $\overline{\varrho}\colon\overline{D}\to \Sigma$ which agrees with $\varrho$ on $S^1$ and has the following property.
For each $i=1,\dots, m$ there exists a $2$--simplex $\sigma_i\subset\Sigma$ such that $\overline{\varrho}$ maps $B_i$ biLipschitz homeomorphically onto $\sigma_i$, and $\overline{\varrho}(\overline{D}\setminus \cup_{i=1}^m\operatorname{int}(B_i))\subset \Sigma^{(1)}$.
 Moreover, if $\varrho|_{S^1}$ is Lipschitz then $\overline{\varrho}$ can be taken to be Lipschitz on $\overline{D}$.
\el

Let $U\subset\R^2$ be open and $f\colon U\to\R^2$ continuous. Let $A\subset\R^2$ be the subset of points $y\in\R^2$ such that $N(f,y)<\infty$. For $y\in A$ and $x\in f^{-1}(y)$ we denote by $\iota(f, x)$ the winding number of the curve $f\circ\gamma$ with respect to $y$, where $\gamma\colon S^1\to \R^2$ is given by $\gamma(z) = x+rz$ and $r>0$ is chosen so small that $\bar{B}(x,r)\subset U$ and $\bar{B}(x,r)\cap f^{-1}(y)=\{x\}$. Clearly, the winding number of $f\circ\gamma$ with respect to $y$ is independent of the choice of such $r$. It follows from \cite[Lemma 5.2]{Rad38} that there exists an at most countable set $N\subset A$ such that $|\iota(f,x)|\leq 1$ for each $y\in A\setminus N$ and every $x\in f^{-1}(y)$.

\begin{proof}[Proof of Lemma~\ref{lem:approx-cont-maps-multiplicity}]
 Denote by $\sigma_1,\dots, \sigma_n$ the finitely many $2$--simplices of $\Sigma$ and notice that $\hm^2(\sigma_i) = |\sigma_i| = \frac{\sqrt{3}}{4}r^2$. By the remark after the statement of the lemma, for each $i=1,\dots, n$ there exists  $y_i\in\interior(\sigma_i)$ such that $$|\sigma_i|\cdot N(\varrho, y_i) \leq \int_{\sigma_i} N(\varrho, y)\,d\hm^2(y)$$ and $|\iota(\varrho, x)|\leq 1$ for every $x\in\varrho^{-1}(y_i)$. Write $\varrho^{-1}(y_i) = \{x_{i,1},\dots, x_{i,m_i}\}$, where $m_i= N(\varrho,y_i)$, and choose $r_i>0$ so small that the balls $\bar{B}(x_{i,j}, 2r_i)$ are contained in $\varrho^{-1}(\interior(\sigma_i))$ and are pairwise disjoint.
 Let $\pi\colon \Sigma\setminus \{y_1,\dots, y_n\}\to \Sigma^{(1)}$ be the continuous map which is the identity on $\Sigma^{(1)}$ and such that $\pi|_{\sigma_i\setminus\{y_i\}}$ is the radial projection onto $\partial \sigma_i$ with projection center $y_i$. 
 
 We let $\overline{\varrho}\colon\overline{D}\to \Sigma$ be the continuous map which agrees with $\pi\circ\varrho$ on the complement of the balls $B(x_{i,j}, 2r_i)$ and such that $\overline{\varrho}|_{\bar{B}(x_{i,j}, 2r_i)}$ is defined as follows. If $\iota(\varrho, x_{i,j})=0$ then $\overline{\varrho}|_{\partial B(x_{i,j}, 2r_i)}$ is contractible in $\partial \sigma_i$ and it therefore has a continuous extension to $\bar{B}(x_{i,j}, 2r_i)$ with image inside $\partial\sigma_i$. If $\iota(\varrho, x_{i,j})=\pm 1$ then $\overline{\varrho}|_{\partial B(x_{i,j}, 2r_i)}$ is homotopic inside $\partial \sigma_i$ to a biLipschitz parametrization of $\partial \sigma_i$. We define $\overline{\varrho}|_{\bar{B}(x_{i,j}, 2r_i)\setminus B(x_{i,j},r_i)}$ to be such a homotopy and we let $\overline{\varrho}|_{\bar{B}(x_{i,j}, r_i)}$ be a biLipschitz homeomorphism onto $\sigma_i$ which extends $\overline{\varrho}|_{\partial B(x_{i,j}, r_i)}$. It is clear that $\overline{\varrho}$ has all the properties listed in the statement of the proposition, except the last one.
 
 In order to prove the last statement, suppose $\varrho|_{S^1}$ is Lipschitz continuous. Set $\Omega:= D\setminus\cup_{i,j} B(x_{i,j},r_i)$ and notice that $\overline{\varrho}(\overline{\Omega})\subset\Sigma^{(1)}$ and the restriction of $\overline{\varrho}$ to $\partial \Omega$ is Lipschitz. Since $\Sigma^{(1)}$ is locally Lipschitz $k$--connected for every $k\geq 0$ it follows from the Lipschitz extension results \cite{LS05} and \cite{Hoh93} that we can approximate $\overline{\varrho}|_{\overline{\Omega}}$ arbitrarily closely by a Lipschitz map with image in $\Sigma^{(1)}$ which agrees with $\overline{\varrho}$ on $\partial \Omega$. This concludes the proof.
\end{proof}

\begin{proof}[Proof of Proposition~\ref{prop:Lipschitz-approx-neighborhood}]
 Let $r>0$ be sufficiently small and let $\Sigma$, $\psi$, $\varphi$ and the constant $C$ be as in Lemma~\ref{lem:approximation-simplicial-complex}. Let $c\colon S^1\to X$ be a constant speed parametrization of $\partial X$ and set $L:= \length(\partial X)$.
 We first claim that there exist a constant $C_1$ only depending on $C$ and a continuous map $\varrho\colon\overline{D}\to \Sigma$ such that $$\int_{\Sigma} N(\varrho,y)\,d\hm^2(y)\leq C_1(\hm^2(X) + Lr)$$ and $\varrho|_{S^1}$ is $C_1L$--Lipschitz with image in $\Sigma^{(1)}$ and $d(\varrho(t), \psi(c(t)))\leq C_1r$ for all $t\in S^1$. Indeed, the $CL$--Lipschitz curve $\psi\circ c$ is homotopic to a $C'L$--Lipschitz curve $\gamma\colon S^1\to \Sigma$ satisfying $\gamma(S^1)\subset \Sigma^{(1)}$ and $$d(\psi(c(t)), \gamma(t))\leq C'r$$ for all $t\in S^1$ via a Lipschitz homotopy $h$ of area $\Area(h)\leq C'Lr$,  where $C'$ only depends on $C$. Such $\gamma$ can be obtained by replacing $\psi\circ c$ on the closure of each component of $(\psi\circ c)^{-1}(\interior(\sigma))$ by the constant speed shortest curve in $\partial\sigma$ for every $2$--simplex $\sigma$. The homotopy $h$ is the straight line homotopy in $\sigma$. In particular, we obtain from the area formula that $$\int_\Sigma N(h,y)\,d\hm^2(y)= \Area(h) \leq C'Lr.$$ By the Jordan-Schoenflies theorem there exists a homeomorphism $\eta\colon \overline{D}\to X$ which extends $c$. The coarea inequality for Lipschitz maps \cite[Theorem 2.10.25]{Fed69} implies $$\int_\Sigma N(\psi\circ\eta,y)\,d\hm^2(y) = \int_\Sigma N(\psi,y)\,d\hm^2(y)\leq C^2\hm^2(X).$$ The map $\varrho$ given by $\varrho(z) = \psi(\eta(2z))$ when $|z|\leq \frac{1}{2}$ and by $h(z/|z|, 2|z|-1)$ when $\frac{1}{2}\leq |z|\leq 1$ satisfies the claim above for some $C_1$ only depending on $C$. 
 
 Next, let $\overline{\varrho}\colon \overline{D}\to \Sigma$ be a Lipschitz map as in Lemma~\ref{lem:approx-cont-maps-multiplicity} associated with the map $\varrho$. We then have $$\Area(\overline{\varrho})\leq \int_\Sigma N(\varrho,y)\,d\hm^2(y)\leq C_1(\hm^2(X) + Lr).$$ Moreover, since $\alpha:=\varphi\circ\overline{\varrho}|_{S^1}$ is $(CC_1L)$--Lipschitz and satisfies 
 \begin{equation*}
      d(\alpha(t),c(t))\leq d(\alpha(t), \varphi(\psi(c(t)))) + d(\varphi(\psi(c(t))), c(t))
      \leq C(C_1 + 1)r
 \end{equation*}
 for all $t\in S^1$, the Lipschitz extension property of $E(X)$ implies that there exists a Lipschitz homotopy $g$ from $\alpha$ to $c$ in $E(X)$ of area bounded by $C''Lr$ and with image in the $(C''r)$--neighborhood of $c$ for some $C''$ only depending on $C$. The Lipschitz map $v\colon \overline{D}\to E(X)$ given by $v(z) = \varphi(\overline{\varrho}(2z))$ when $|z|\leq \frac{1}{2}$ and by $v(z) = g(z/|z|, 2|z|-1)$ when $\frac{1}{2}\leq |z|\leq 1$ agrees with $c$ on $S^1$; its image is contained in the $C_2r$--neighborhood of $X$, and $$\Area(v)\leq \Area(\varphi\circ\overline{\varrho}) + \Area(g)\leq C^2C_1(\hm^2(X) + Lr) + C''Lr \leq C_2(\hm^2(X) + Lr)$$ for some constant $C_2$ only depending on $C$. The proposition now follows.
\end{proof}

\begin{proof}[Proof of Corollary~\ref{cor:intro-rectifability}]
 Let $\Omega\subset X$ be a Jordan domain with finite boundary length. Such $\Omega$ can be constructed as in the proof of Lemma~\ref{lem:biLip-Jordan}. By Theorem~\ref{thm:non-trivial-Sobolev-Jordan-intro} and its proof there exists $u\in \Lambda(\partial \Omega, \overline{\Omega})$. We claim that $\Area(u)>0$. Indeed, otherwise the infimum of energies over all maps in $\Lambda(\partial\Omega, \overline{\Omega})$ would be zero by Morrey's $\varepsilon$--conformality Lemma \cite{FW20}. Hence an energy minimizer, which exists by Theorem~\ref{thm:existence-energy-min}, would have zero energy and would thus be constant, a contradiction. This proves the claim.
 
Let $A_1\subset A_2\subset\dots\subset D$ be measurable sets with $|D\setminus\bigcup A_i|=0$ and such that the restriction $u|_{A_i}$ is Lipschitz for every $i$, see e.g.~\cite[Proposition 3.2]{LW15-Plateau}. Since $\Area(u)>0$ there exists $i$ such that $\Area(u|_{A_i})>0$. By the area formula we have $$\Area(u|_{A_i}) = \int_{u(A_i)}N(u|_{A_i}, x)\,d\hm^2(x)$$ and hence $\hm^2(u(A_i))>0$. This completes the proof.
\end{proof}

\section{Finishing the proofs of the almost parametrization results}\label{sec:proofs-almost-param}

In this section we finish the proofs of the almost parametrization results given in the introduction and discuss some additional consequences.

\begin{proof}[Proof of Theorem~\ref{thm:main-almost-param}]
 Denote by $d$ the metric on $X$ and consider the length metric $d_{\overline{\Omega}}$ on $\overline{\Omega}$. It follows as in the second part of the proof of Theorem~\ref{thm:non-trivial-Sobolev-Jordan-intro} that the metric space $Y=(\overline{\Omega}, d_{\overline{\Omega}})$ is geodesic, homeomorphic to $\overline{D}$, and $\length(\partial Y)$ and $\hm^2(Y)$ are finite. Hence, $\Lambda(\partial Y, Y)$ is not empty by Theorem~\ref{thm:non-trivial-Sobolev-Jordan-intro}. Therefore, by Theorem~\ref{thm:existence-energy-min}, there exists an energy minimizer $v$ in $\Lambda(\partial Y, Y)$ and every such $v$ is infinitesimally $\frac{4}{\pi}$--quasiconformal. Theorem~\ref{thm:cont-energy-min-intro} further implies that $v$ has a continuous representative which continuously extends to the boundary, denoted by $v$ again. Finally, Theorem~\ref{thm:energy-min-monotone} implies that $v$ is monotone.
 
 Since the identity $\pi\colon Y\to (\overline{\Omega}, d)$ is a homeomorphism the map $u\colon \overline{D}\to \overline{\Omega}\subset X$ defined by $u\coloneqq \pi\circ v$ is continuous, surjective, and monotone. Since $\pi$ is $1$--Lipschitz and its restriction to $Y\setminus \partial Y$ is a local isometry it follows furthermore that $u$ belongs to $N^{1,2}(D, X)$ and that $u$ is infinitesimally $\frac{4}{\pi}$--quasiconformal. Thus, $u$ satisfies \eqref{eq:analytic-qc-mod} by the discussion at the beginning of Section~\ref{sec:properties-qc-almost-parametrizations}. This completes the proof.
\end{proof}

\br\label{rem:loc-geod-weakened}
 {\rm The condition in Theorem~\ref{thm:main-almost-param} that the metric space $X$ be locally geodesic can be relaxed. It suffices to assume that $X=(X,d)$ is rectifiably connected, the length metric $d_i$ induces the same topology, and $X_i=(X,d_i)$ has locally finite Hausdorff $2$--measure. To see that this suffices, we first observe that the family $\Lambda$ of maps $u\in N^{1,2}(D, X)$ such that $u$ is the uniform limit of homeomorphisms $\overline{D}\to \overline{\Omega}$ is not empty. Indeed, by Theorem~\ref{thm:main-almost-param}, there exists a quasiconformal almost parametrization $v$ of $(\overline{\Omega}, d_i)$. Let $\pi\colon (\overline{\Omega}, d_i)\to (\overline{\Omega}, d)$ be the identity map and notice that the map $u:=\pi\circ v$ is the uniform limit of homeomorphisms $\overline{D}\to\overline{\Omega}$. Since $\pi$ is $1$--Lipschitz it follows that $u$ belongs to $N^{1,2}(D,\overline{\Omega})$, so $\Lambda$ is not empty. Next, one shows that $\Lambda$ contains an energy minimizer. For this, let $(u_n)\subset\Lambda$ be an energy minimizing sequence. After precomposing with M\"obius transformations we may assume that the $u_n$ satisfy a $3$--point condition and so, by \cite[Proposition 7.4]{LW15-Plateau}, the sequence $(u_n|_{S^1})$ is equi-continuous.  Hence, the proof of Proposition~\ref{prop:essential-oscillation-energy-min} shows that the sequence $(u_n)$ is equi-continuous. Thus, after passing to a subsequence, we may assume that $(u_n)$ converges uniformly to a map $u$. This map belongs to $\Lambda$ and is an energy minimizer in $\Lambda$ and thus is infinitesimally $\frac{4}{\pi}$--quasiconformal by \cite{LW16-energy-area}. In particular, $u$ satisfies \eqref{eq:analytic-qc-mod}, see Section~\ref{sec:properties-qc-almost-parametrizations}.}
\er

The following variant of Theorem~\ref{thm:main-almost-param} is closer to the statement of the Riemann mapping theorem.

\bc\label{cor:almost-param-Riem-map-version}
 Let $X$ be a locally geodesic metric space homeomorphic to $\R^2$ and of locally finite Hausdorff $2$--measure. If $U\subset X$ is an open and simply connected set with compact closure then there exists a continuous, monotone surjection $u\colon D\to U$ such that 
 \begin{equation}\label{eq:analytic-qc-mod-cor-Riemann-map}
  \MOD(\Gamma) \leq K\cdot\MOD(u\circ \Gamma)
 \end{equation}
 for every family $\Gamma$ of curves in $D$, where $K=\frac{4}{\pi}$.
\ec

\begin{proof}
 Let $\Omega\subset X$ be a Jordan domain containing $\overline{U}$. We can approximate $\partial \Omega$ by a biLipschitz Jordan curve as in the proof of Lemma~\ref{lem:biLip-Jordan} and may therefore assume that $\Omega$ is of finite boundary length. By Theorem~\ref{thm:main-almost-param}, there exists a continuous, surjective, monotone map $v\colon\overline{D}\to\overline{\Omega}$ such that $\MOD(\Gamma)\leq \frac{4}{\pi}\cdot\MOD(v\circ\Gamma)$ for every family $\Gamma$ of curves in $\overline{D}$. Since $v$ is monotone it follows that $V:=v^{-1}(U)\subset D$ is also simply connected, see \cite[Section 2.3]{LW20-param}. By the Riemann mapping theorem, there exists a conformal diffeomorphism $\varphi\colon D\to V$. Then the map $u\colon D\to U$ given by $u:=v\circ\varphi$ has the desired properties.
\end{proof}

\begin{proof}[Proof of Corollary~\ref{cor:geom-qc-param-reciprocal}]
 Let $U\subset X$ be a Jordan domain containing $\overline{\Omega}$. Arguing as in the proof of Corollary~\ref{cor:almost-param-Riem-map-version}, we may assume that $U$ has finite boundary length. By Theorem~\ref{thm:main-almost-param}, there exists a continuous, surjective, monotone map $v\colon\overline{D}\to\overline{U}$ such that $\MOD(\Gamma)\leq \frac{4}{\pi}\cdot\MOD(v\circ\Gamma)$ for every family $\Gamma$ of curves in $\overline{D}$. It follows from  Propositions~\ref{prop:points-modulus-zero-implies-homeo} and \ref{prop:reciprocal-implies-second-modulus-inequality} that $v$ is geometrically quasiconformal homeomorphism. Finally, by the Riemann mapping theorem, there exists a conformal diffeomorphism $D\to \Omega'$, where $\Omega' = v^{-1}(\Omega)$, which moreover extends to a homeomorphism $\varrho\colon\overline{D}\to \overline{\Omega'}$. The composition $u:= v\circ\varrho\colon\overline{D}\to\overline{\Omega}$ is a geometrically quasiconformal homeomorphism.
\end{proof}

Another consequence of Theorem~\ref{thm:main-almost-param} is the following variant for discs of the Bonk-Kleiner quasisymmetric uniformization theorem \cite{BK02}, see also \cite{LW20-param}.

\bc\label{cor:qs-param}
 Let $X$ be a geodesic metric space homeomorphic to $\overline{D}$ and with finite boundary length. If there exists $L>0$ such that $$\hm^2(B(x,r))\leq Lr^2$$ for all $x\in X$ and $r>0$ and if $X$ is linearly locally connected then there exists a quasisymmetric homeomorphism $u\colon\overline{D}\to X$.
\ec

Using a quasisymmetric gluing theorem exactly as in the proof of \cite[Proposition 6.4]{LW20-param} one obtains an analogous statement when $X$ is homeomorphic to $S^2$ and thus the Bonk-Kleiner quasisymmetric uniformization theorem \cite{BK02}.

\begin{proof}
 By Theorem~\ref{thm:main-almost-param}, there exists a continuous, monotone surjection $u\colon\overline{D}\to X$ satisfying \eqref{eq:analytic-qc-mod}. Since the quadratic upper bound for the Hausdorff measure of balls implies \eqref{eq:points-zero-modulus} by \cite[Lemma 7.18]{Hei01} it follows from Proposition~\ref{prop:points-modulus-zero-implies-homeo} that $u$ is a homeomorphism. Finally, $u$ is quasisymmetric by Proposition~\ref{prop:upgrade-qs}.
\end{proof}

\def\cprime{$'$} \def\cprime{$'$} \def\cprime{$'$} \def\cprime{$'$}

\end{document}